\documentclass[a4paper,12pt,leqno]{amsart}

\usepackage[latin1]{inputenc}
\usepackage{utopia}
\usepackage{mathptm}

\usepackage{enumerate}
\usepackage{latexmac}
\usepackage{mes_abbrev} 

\newtheorem{theorem}{Theorem}

\newtheorem{proposition}[theorem]{Proposition}

\theoremstyle{definition}

\theoremstyle{remark}

\newcommand{\gesp}[1]{\mathbf{E}^{(g)}\etc{#1}}

\newcommand{\gPP}{\mathbf{P}^{(g)}}

\newcommand{\gespcro}[1]{\gesp{\crochet{#1}}}

\newcommand{\qud}{q_{12}}
\newcommand{\qut}{q_{13}}

\newcommand{\tz}{\tilde{Z}}
\newcommand{\ta}{\tilde{\alpha}}

\newcommand{\fpar}[2]{\frac{\partial #1}{\partial #2}}

\newcommand{\lin}[1]{{#1}^{\mathbf{(lin)}}}
\newcommand{\zlin}{\lin{Z}}
\newcommand{\qlin}{\lin{q}}
\newcommand{\alphalin}{\lin{\alpha}}
\newcommand{\hpsx}{h+\sqrt{x} }
\newcommand{\hpsxg}{\hpsx g}

\newcommand{\tbz}{\bar{Z}}
\newcommand{\ba}{\bar{\alpha}}
\newcommand{\SK}{Sherrington-Kirkpatrick\xspace}
\newcommand{\SKs}{Sherrington-Kirkpatrick's\xspace}

\begin{document}
\title[High temperature Sherrington-Kirkpatrick model for general spins]{High temperature Sherrington-Kirkpatrick model for general spins}

\author{Philippe Carmona}
\address{Philippe Carmona\\
Laboratoire Jean Leray, UMR 6629,
Universit{\'e} de Nantes, BP 92208\\
F-44322 Nantes Cedex 03
BP 
}
\email{philippe.carmona@math.univ-nantes.fr}
\begin{abstract}
  Francesco Guerra and Fabio Toninelli~\cite{guerra02:_quadr_sherr_kirkp,guerra2000:sum_rules} have developped a very powerful technique to study the high temperature behaviour of the \SK mean field spin glass model. They show that this model is asymptoticaly comparable to a linear model. The key ingredient is a clever interpolation technique between the two different Hamiltonians describing the models. 

This paper contribution to the subject are the following:
\begin{itemize}
\item The replica-symmetric solution holds for general spins, not just $\pm 1$ valued.
\item The proof does not involve cavitation but only first order differential calculus and Gaussian integration by parts.
\end{itemize}

\end{abstract}



\date{\today}\vfuzz=4pt \hfuzz=40pt
\maketitle

\section{Introduction}

In spin glasses models a generic configuration $\sigma=(\sigma_1,\ldots,\sigma_n)$ represents the position of $n$ spins $\sigma_=\pm 1$. In the \SK model, the external disorder is given by $n(n-1)/2$ iid (independent identically distributed) random variables $(g_{ij})_{1\le i<j\le n}$ assumed $\Nrond(0,1)$ that is centered unit Gaussian. The Hamiltonian, for a given inverse temperature $\beta$ and in some external field of strength $h$, is given by
$$H^{(SK)}_n(\sigma) =H^{(SK)}_n(\sigma,g) = \beta\unsur{\sqrt{n}} \sum_{i<j} g_{ij} \sigma_i\sigma_j + h \sum_{i=1}^n \sigma_i\,.$$
 
 The partition function $Z^{(SK)}_n(\beta)$ and free energy $\alpha^{(SK)}_n(\beta)$ of the model are given by
$$ Z^{(SK)}_n(\beta) = \sum_\sigma e^{ H_n(\sigma)}\,,\quad\alpha_n^{(SK)}(\beta) = \unsur{n} \gesp{\log Z_n(\beta)}\,.$$

Physicists~\cite{mezard:_spin} and mathematicians~\cite{ talsaintflour2000} both proved that for high temperatures, $\alpha_n^{(SK)}(\beta)$ converges to the \SK replica symmetric solution
$$ \alpha_n^{(SK)}(\beta) \to \alpha_\infty^{(SK)}(\beta) = \log 2 + \frac{\beta^2}{4} (1-q)^2 + \gesp{\log \cosh( h + g \beta {\sqrt{q}})}\,,$$

where $q$ is the unique solution of the equation $q=\gesp{\tanh^2(\beta g \sqrt{q} + \beta h)}$ (here and in the following $g$ is a standard $\Nrond(0,1)$ Gaussian random variable).

To introduce the concept of general spins, we need to normalize things
\begin{itemize}
\item $\beta=\sqrt{2t}$.
\item We divide the partition function by the number $2^n$ of configurations.
\item We compensate each weight (Boltzmann factor) $e^{ H^{(SK)}_n(\sigma)}$ so that it has expectation $1$ with respect to the external disorder.
\end{itemize}
With the notation $H_n(\sigma)= \sqrt{\frac{2}{n}} \sum_{i<j} g_{ij}\,\sigma_i\sigma_j + \unsur{\sqrt{n}} \sum_i g_{ii} \sigma_i^2$ the partition function is 

$$Z_n(t) = \esp{e^{ \sqrt{t}H_n(\sigma) + h \sum_{i=1}^n \sigma_i - \frac{t}{2}n}} = 2^{-n} \sum_{\sigma} e^{ \sqrt{t}H_n(\sigma) + h \sum_{i=1}^n \sigma_i - \frac{t}{2}n}\,,$$

where under the probability $\PP$, $\sigma_i$ are iid random variables with distribution $\prob{\sigma_i=\pm 1} = \undemi$. We have now
$$ \alpha_n(t) = \unsur{n}\gesp{\log Z_n(t)} = \alpha_n^{(SK)}(\beta) -\log 2 -\frac{t}{2}\,,$$

and thus \SKs result can be rephrased

$$ \a_n(t) \to \a_\infty(t) = \frac{t}{2} q^2 + \gesp{\log \cosh (h + g\sqrt{2qt})} - tq\,,$$ with $q$ the unique solution of $q=\gesp{\tanh^2(\beta g \sqrt{q} + h)}$.

(let us observe that the introduction of a fixed random variable $\sum_i g_{ii} \sigma_i^2=\sum_i g_{ii}$ does not change the free energy and simplifies the computations). 

\smallskip

To generalize the model we
assume now that the spins are not $\pm 1$ valued, but that they are just, under $\PP$ symmetric iid random variables with values in $\etc{-1,1}$ ; for instance uniformly ditributed on $\etc{-1,1}$. We introduce, as usual, the mutual overlap between two spins $\sigma,\tau$

$$ \qud=\qud(\sigma,\tau) = \unsur{n} \sum_{i=1}^n \sigma_i \tau_i\,.$$
Then 
$$H_n(\sigma)=\sqrt{\frac{2}{n}} \sum_{i<j} g_{ij}\,\sigma_i\sigma_j+\unsur{\sqrt{n}} \sum_{i=1} g_{ii} \sigma_i^2\,,$$ is a centered Gaussian process with covariance
\begin{align*}
 \esp{H_n(\sigma)H_n(\tau)} &= \frac{2}{n} \sum_{i<j} \sigma_i \sigma_j \tau_i \tau_j + \unsur{n} \sum_i \sigma_i^2 {\tau_i}^2 \\
&= n \qud^2(\sigma,\tau)\,.
\end{align*}

Accordingly,

$$ Z_n(t) = \esp{e^{\sqrt{t} H_n(\sigma) + h \sum_{i=1}^n \sigma_i - \frac{t}{2} n \qud^2(\sigma,\sigma)}}\,,\quad \alpha_n(t) = \unsur{n}\gesp{\log Z_n(t)}\,.$$

Our main result is the following

\begin{theorem}
  Let $\phi(u,v)$ and $\qlin(x)$ be defined by
$$ e^{\phi(u,v)} = \esp{e^{ u\sigma_i + v \sigma_i^2}}\,,\quad \qlin(x) =\gesp{\partial_u\phi^2(h+g\sqrt{x}, -x/2)}\,.$$

There exists a number $t_c>0$ such that for all $t\le t_c$
$$\a_n(t) \to\a_\infty(t) = \frac{t}{2} q^2 + \gesp{\phi(h + g\sqrt{2qt},-qt)}\,,$$

where $q=q_c(t)$ is the unique solution of $\qlin(2qt)=q$.
\end{theorem}

For the classical \SK model, $\phi(u,v) = v+\log \cosh u$ and $\partial_u\phi^2(u,v) = \tanh^2(u)$.

\medskip

Let us now explain what is the key ingredient of the proof : Francesco Guerra's interpolation technique. It has been successfully used by Guerra~\cite{guerra2000:sum_rules}, Guerra and Toninelli~\cite{guerra02:_quadr_sherr_kirkp} and Talagrand~\cite{talagran02:_high_temper_sherr_kirkp}, to show that the replica-symmetric formula holds in a region which probably coincides with the Almeida-Thouless region.

We introduce a simpler model, with a linear random Hamiltonian $\sqrt{x} \Lambda_n(\sigma) + h \sum_i \sigma_i$, with $\Lambda_n(\sigma)=\sum_i J_i \sigma_i$, and $J_i$ are independent $\Nrond(0,1)$. For this linear model, it is immediate to compute the free energy $\alphalin(x)$. 

Then, we consider a two parameter Hamiltonian, $\sqrt{t} H_n(\sigma)+\sqrt{x} \Lambda_n(\sigma) + h \sum_i \sigma_i$ and compare the free energies obtained for $x=0$, the \SK model, and $t=0$ the linear model. This is easily done for some $x=2qt$ with $q$ a solution of an equation involving $t$ and $h$.

\section{The linear model}
\label{sec:linear-model}
This model is simpler to study than the \SK model, because here the partition function is a product of independent factors.

\medskip

Recall that the spins $\sigma_i$ are assumed to be independent identically distributed, with values in $\etc{-1,1}$. The function $\phi$ denotes the mixed Laplace exponent
$$ \phi(u,v) = \log \esp{e^{u\sigma_i + v \sigma_i^2}}\,.$$


The Hamiltonian is a linear random form on the spins $\Lambda_n(\sigma)=\sum_{i=1}^n \sigma_i J_i$ with $(J_i)_i$ independent $\Nrond(0,1)$ : its covariance is $n$ times the the overlap:
$$\esp{\Lambda_n(\sigma)\Lambda_n(\sigma')}=n\qud(\sigma,\sigma')\,.$$ 
There is an external linear field of strength $h\ge 0$ so the partition function is:

$$\zlin_n(x) = \esp{e^{\sqrt{x}\sum_{i=1}^n \sigma_i J_i -\frac{x}{2}n\qud(\sigma,\sigma)+ h \sum_{i=1}^n \sigma_i }}\,.$$

\begin{proposition}
\label{sec:linear-model:mean-overlap}
 1) The mean overlap is given by
  \begin{equation}
    \label{eq:linear:formeqlin}
    \gesp{\crochet{\qud}}=\gesp{\partial_u\phi^2(h+\sqrt{x} g,-x/2)}\egaldef\qlin(x)\,.
  \end{equation}
2) There exists two numbers $\l_0>0$ and $L$ such taht for any $\l \le \l_0$, and any $x,h,n$:
\begin{equation}
  \label{eq:linear:exponential}
\gesp{\log\crochet{\exp \l n\etp{\qud(\sigma,\tau) -\qlin(x)}^2}} \le L\,.
\end{equation}
where the bracket represents a double integral with respect to the Gibbs measure.
\end{proposition}

\begin{proof}
1)  Let $\alphalin(x) = \unsur{n} \gesp{\log Z_n(x)}$. On the one hand, integration by parts (see section~\ref{sec:ipp}) yields

$$  \frac{d \alphalin_n}{dx} = -\undemi\gespcro{\qud}.$$

On the other hand, a direct computation using the independence of spins yields
$$ \zlin_n(x) = \prod_{i=1}^n \esp{e^{(h+\sqrt{x} J_i)\sigma_i -\frac{x}{2}\sigma_i^2}} = \exp \sum_{i=1}^ n \phi(h+\sqrt{x} J_i,-\frac{x}{2})\,.$$ 

Therefore,
$ \alphalin(x) = \gesp{\phi(h+\sqrt{x} g,-\frac{x}{2})}$, and by integration by parts 

\begin{align*}
\frac{d \alphalin_n}{dx} &= \gesp{\unsur{2\sqrt{x}}\partial_u\phi(\hpsxg,-\frac{x}{2})} -\undemi \gesp{\partial_v\phi(\hpsxg,-\frac{x}{2})}\\
&= \gesp{\undemi\partial^2_{u^2}\phi(\hpsxg,-\frac{x}{2})} -\undemi \gesp{\partial_v\phi(\hpsxg,-\frac{x}{2})}
\end{align*}

Hence,
$$ \qlin(x) = \gesp{\crochet{\qud}} = \gesp{ (\partial_v - \partial^2_{u^2})(\hpsxg,-\frac{x}{2})} = \gesp{\partial_u\phi^2(h+\sqrt{x} g,-x/2)}$$
because
$\partial_v\phi - \partial^2_{u^2}\phi=(\partial_u\phi)^2$.

\medskip
2) This part of the proof can be established via the cavitation techniques introduced by Talagrand. This is a direct consequence of a stronger  exponential inequality : there exists a constant $C>0$, such that for all $x,h$ ($h$ small enough)

$$ \gespcro{\exp \frac{n}{C} (\qud(\sigma,\tau) - \qlin(x))^2} \le C\,.$$

However, we shall give a more direct proof, which avoids the hassles of cavitation.

\def\bsig{\mathbf{\sigma}}
\def\btau{\mathbf{\tau}}
\newcommand{\espgamma}[1]{\mathbf{E}^\gamma\etc{#1}}

\begin{gather*}
\text{Let}\quad V_n(\bsig) = \sqrt{x} \sum_{i=1}^n \sigma_i J_i -n \frac{x}{2} \qud(\bsig,\bsig) + h \sum_{i=1}^n \sigma_i\,,\\
\text{and}\quad U_n(\bsig,\btau) = V_n(\bsig) + V_n(\btau) + \l (q-\qud(\bsig,\btau))^2\,.
\end{gather*}

If $Z_n(x,\l)$ denotes the partition function $Z_n(x,\l)=\esp{U_n(\bsig,\btau)}$ then we shall prove that for $0\le \l \le \unsur{20}$ and $q=\qlin(x)$,

$$\gesp{\log Z_n(x,\l) -\log Z_n(x,0)} \le L\,,$$

where $L$ is a number.

\smallskip

To    this end, observe that by introducing an auxiliary unit gaussian random variable $\gamma$ with associated expectation $\mathbf{E}^\gamma$, we can write, using the independence of spins,
\begin{align*}
Z_n(x,\l) &= \espgamma{\esp{ e^{V_n(\bsig) + V_n(\btau) + \sqrt{2\l n} \gamma (\qud(\bsig,\btau) -q)}}} \\
&= \espgamma{ e^{-q \gamma \sqrt{2\l n}} \exp(\sum_{i=1}^n \psi(h+\sqrt{x} J_i,-x/2,\l))}
\end{align*}
with
$$ e^{\psi(u,v,\l)} = \esp{e^{\sigma_1 u + \sigma_1^2 v + \tau_1 u + \tau_1^2 v +\l \sigma_1 \tau_1}}\,.$$

It is easily seen that the function $\l \to \psi(u,v,\l)$ is convex, twice differentiable, of first derivative at $\l=0$
$$\frac{\partial\psi}{\partial\l}(u,v,\l=0)= (\partial_u \phi)^2(u,v) \,,$$
and satisfies
$$ 0 \le \frac{\partial^2 \psi}{\partial^2\l^2} (u,v,\l) \le 4\,.$$
Indeed $\psi(u,v,\l)-\psi(u,v,0) $ is the logarithm of the Laplace transform in $\l$ of a random variable taking its values in $\etc{-1,1}$. Therefore, the second derivative is the variance with respect to a twisted probability measure of a random variable taking its values in $\etc{-1,1}$, and is bounded by $4$.

Therefore,
$$ \psi(u,v,\l) \le \psi(u,v,0) + \l (\partial_u \phi)^2(u,v) + 2 \l^2\,.$$
Hence,

$$ \frac{Z_n(x,\l)}{Z_n(x,0)} \le \espgamma{\exp \etp{ -q \gamma \sqrt{2\l n} + 4 \l \gamma^2 + \sum_{i=1}^n \gamma \sqrt{2\l/n} (\partial_u \phi)^2(h+\sqrt{x}J_i, -x/2)}}\,.$$

From now on, $L$ is a number whose value may change from line to line. One easily shows (e.g. with the help of H\"older's inequality) that for $0\le v\le 1/5$ and any $u$
$$ \espgamma{e^{\gamma u + v \gamma^2}} \le L e^{u^2}\,.$$

this entails that for $0\le \l \le 1/20$,
$$ \log \frac{Z_n(x,\l)}{Z_n(x,0)} \le L + \frac{2\l}{n} \etp{\sum_{i=1}^n (\partial_u \phi)^2(h+\sqrt{x}J_i, -x/2)-q}^2\,.$$

The random variables $X_i= (\partial_u \phi)^2(h+\sqrt{x}J_i, -x/2)-q$ are independent identically distributed. Since $q=\qlin(x)$, they are centered, and we observe that they are bounded : $X_i\in\etc{-1,1}$. We can now conclude

$$ \gesp{ \log \frac{Z_n(x,\l)}{Z_n(x,0)}} \le L +  \frac{2\l}{n} \gesp{(X_1 +\cdots + X_n)^2} \le L + 2\l \le L'\,.$$
\end{proof}
\section{The Sherrington Kirkpatrick model}
\xcom{This is a spin-glass model with configuration $X(n)=(\sigma_1,\ldots,\sigma_n)$ with the $(\sigma_i)_i$ iid. We assume that the spin has exponential moments of all order and we define the Laplace exponent 
$$e^\psi(\l)=\esp{e^{\l \sigma_1}}<+\infty \qquad(\l \in\R)\,.$$
}

The Hamiltonian is a bilinear random form on the spins
$$ H_n(\sigma)=\sqrt{\frac{2}{n}}\sum_{i<j} g_{ij} \sigma_i\sigma_j + \unsur{\sqrt{n}} \sum_{i=1}^n g_{ii} \sigma_i^2\,,$$
where $(g_{ij})_{i,j}$ is a family of independent gaussian $\Nrond(0,1)$. If, as usual, $\qud$ denotes the overlap
$$ \qud=\qud(\sigma,\tau) = \unsur{n}\sum_{i=1}^n \sigma_i \tau_i\,,$$
then the covariance of the process $(H_n(\sigma))_\sigma$ is given by
$$\gesp{H_n(\sigma) H_n(\tau))}=n\qud^2(\sigma,\tau) \,.$$

We want to understand the asymptotic behaviour of the partition function
$$ Z_n(t) = \esp{e^{\sqrt{t} H_n(\sigma) -\frac{t}{2}n\qud^2(\sigma,\sigma)+ h \sum_{i=1}^n \sigma_i}}$$
and so we introduce the mean free energy
$$ \alpha_n(t) = \unsur{n}\gesp{\log Z_n(t)}\,.$$

The main idea, taken from Guerra's papers, is to show that in the presence of a linear external field, the partition function $Z_n(t)$ can be compared asymptotically to the partition function of a linear model (see section~\ref{sec:linear-model}) $\zlin_n(x)$ with $x$ an implicit function of $t$. Let us first explain this by giving an upper bound.

\begin{theorem}
  \label{thm:sk:bornesup}
1) With the notations $\alphalin(x)$ and $\qlin(x)$ introduced in section~\ref{sec:linear-model}, we have the upper bound:
$$\limsup_{n\to +\infty} \alpha_n(t) \le \inf_{q>0} \alphalin(2qt)+\frac{t}{2}q^2\,.$$

2) There exists $t_c>0$  such that for $t\le t_c$ and $h>0$, the infimum is attained at $q_c(t)$ the unique solution of $\qlin(2qt)=q$.
\end{theorem}

\begin{proof}
  Following Guerra, we introduce an interpolation between the two Hamiltonians $H_n$ and $\Lambda_n$.The two Gaussian processes $(H_n(\sigma))_\sigma$ and $(\Lambda_n(\sigma))_\sigma$ are assumed independent in the partition function

$$\tz_n(t,x) = \esp{e^{\sqrt{t}H_n + \sqrt{x} \Lambda_n -\frac{n}{2}(x \qud(\sigma,\sigma) + t \qud^2(\sigma,\sigma)) + h\sum_{i=1}^n \sigma_i}}$$
According to gaussian integration by parts (see section~\ref{sec:ipp}), the partial derivatives of 
$\ta_n(t,x)= \unsur{n}\gesp{\log \tz_n(t,x)}$ are given by:

$$ \fpar{\ta_n}{t} = -\undemi\gespcro{\qud^2}\,,$$
and ( as we have already computed in section~\ref{sec:linear-model})
$$\fpar{\ta_n}{x} =-\undemi\gespcro{\qud}\,.$$

Therefore, if we fix $t,x_0,q>0$ and move along the trajectory $x(s) = x_0 - 2 q s$, we have
\begin{equation}
  \label{eq:sk:calculderivee}
\frac{d}{ds} \ta_n(s,x(s)) = \undemi\etp{q^2 - \gespcro{(\qud -q)^2}}
\le \undemi q^2\,.
\end{equation}

Hence,
$$\ta_n(t,x(t)) \le \ta_n(0,x_0) + \frac{t}{2}q^2 = \alphalin(x_0)+\frac{t}{2}q^2\,.$$
We now impose the relationship $x_0=2qt$, in order to get $x(t)=0$ and obtain the upper bound:
$$ \alpha_n(t)=\ta_n(t,0) \le  \alphalin(2qt)+\frac{t}{2}q^2\,.$$
We obtain the desired result by taking $\limsup$ and optimizing in $q>0$.

\medskip
2)
  It is just a matter of computing the derivative of $f(q) = \alphalin(2qt)+\frac{t}{2}q^2$. Using the computations of section~\ref{sec:linear-model}, we get:
$$ f'(q) = t(q - \qlin(2qt))\,.$$

Fact 1 : since $h>0$ , and $\sigma_i$ is not identicaly $0$, we have $\qlin(0)>0$.

Indeed, $\partial_u \phi(u,0)= \frac{\esp{\sigma_1 e^{h \sigma_1}}}{\esp{e^{h \sigma_1}}}$ and therefore, by symmetry

$$ \esp{\sigma_1 e^{h \sigma_1}} = 2 \esp{\sigma_1 \sinh(h \sigma_1)} >0\,.$$

Fact 2 : if $q$ is large enough $f'(q) >0$.
Indeed, the function $\qlin$ is bounded $\valabs{\qlin(x)}\le 1$ since the spins are themselves bounded by $1$.

We have $f'(0)<0$ and $f'(q)>0$ if $q$ is large enough, therefore, all we have to do now is to prove that the equation $q=\qlin(2qt)$ has a unique solution (the infimum of $f$ will then be attained there).

We compute the derivative
$$ \frac{d\qlin(x)}{dx} = -\undemi \gesp{(\partial^4_{u^4} + \partial^2_{uv})\phi(h+g\sqrt{x},-x/2)}\,.$$
Since $\phi$ and all its partial derivatives are bounded, the derivative of $\qlin$ is itself bounded by a constant $C$ :  $\valabs{\frac{d\qlin(x)}{dx}}\le C$. Hence,  for $t<\unsur{4C}$ the function $q\to \qlin(2qt)$ is Lipschitz with constant at most $\undemi$, and therefore it has a unique fixed point.

\end{proof}

To show that this upper bound yields asymptotically a lower bound, the proof is a little more involved.

\begin{theorem}
  \label{thm:sk:inf}
There exists $t_c>0$ such that for $t<t_c$  we have the convergence
$$\lim_{n\to +\infty} \alpha_n(t) = \alphalin(2qt)+\frac{t}{2}q^2  \,,$$ 
where $q=q_c(t)$ is the unique solution of $\qlin(2qt)=q$.
\end{theorem}

\begin{proof}
  We consider now two independent copies $\sigma,\tau$ of the spins, and  the partition function

\begin{gather*}
\tbz(t,x,\l)=\esp{e^{U_n(\sigma,\tau)}}\,\\
U_n(\sigma,\tau)= V_n(\sigma) + V_n(\tau) + \l n (q-\qud(\sigma,\tau))^2 \\
V_n(\sigma)= \sqrt{t} H_n(\sigma) + \sqrt{x} \Lambda_n(\sigma) -\frac{n}{2}(x \qud(\sigma,\sigma) + t \qud^2(\sigma,\sigma))\,.
\end{gather*}

\xcom{
$$\tbz(t,x,\l) = \esp{\exp\etp{
\begin{array}{c}
\sqrt{t}(H_n(\sigma^1)+H_n(\sigma^2)) -\frac{n}{2} t(\qud^2(\sigma^1,\sigma^1) +\qud^2(\sigma^2,\sigma^2))\\
+ \sqrt{x} (\Lambda_n(\sigma^1)+\Lambda_n(\sigma^2)) 
-\frac{n}{2} x(\qud(\sigma^1,\sigma^1) +\qud(\sigma^2,\sigma^2))\\
+{ h \sum_{i=1}^n \sigma_i + \l (q-\qud(\sigma^1,\sigma^2))^2}
\end{array}
}}
,$$
}
where $\l$ is a new positive parameter, and $q$ will be specified later. Accordingly, if $ \ba_n(t,x,\l) = \unsur{2n} \gesp{\log \tbz(t,x,\l)} $ then
\begin{gather*}
  \fpar{\ba_n}{t} = \undemi\gespcro{  \qud^2(\sigma,\tau) - 2\qud^2(\sigma,\sigma')}=\undemi\gespcro{\qud^2 - 2 \qut^2} \\
 \fpar{\ba_n}{x}=\undemi\gespcro{\qud - 2 \qut}\\
\fpar{\ba_n}{\l}= \undemi\gespcro{ (q-\qud)^2}\,.
\end{gather*}
where we have now four replicas $\sigma,\tau,\sigma',\tau'$ with $\sigma,\tau$ coupled by $\l (q-\qud(\sigma,\tau))^2$ and $\sigma',\tau'$ coupled by $\l (q-\qud(\sigma',\tau'))^2$, and $\qut=\qud(\sigma,\sigma')$.

\subsubsection*{Fact 1} For every $\l\ge 0,t\ge 0,x_0\ge 2qt,q\ge 0$ we have:
\begin{equation}
  \label{eq:faitun}
  \ba_n(t,x_0-2qt,\l) \le \ba_n(0,x_0,\l+t) + \frac{t}{2}q^2\,.
\end{equation}
Fix $x_0,t,\l_0>0$. We shall compute the derivative along the trajectory $x(s) = x_0-2 q s$, $\l =\l_0-s$:
$$ \frac{d}{ds} \ba_n(s,x(s),\l(s)) = \undemi q^2 -\gespcro{(q-\qut)^2}
\le \undemi q^2\,.$$
Therefore, integrating between $0$ and $t$,

$$\ba_n(t,x(t),\l(t)) - \ba_n(0,x_0,\l_0) \le \frac{t}{2} q^2\,,$$
and this is the desired inequality.

\xcom{
\subsubsection*{Fact 2} For every $0\le x_0\le x_c,h\le h_c$, $\l,t$ such that $\l \le 1/L$ if $q=\qlin(x_0)$ then,
\begin{equation}
  \label{eq:faitdeux}
  \ba_n(0,x_0,\l) \le \ba_n(0,x_0,0) + \unsur{2n} \log L\,.
\end{equation}

Indeed, by Proposition~\ref{sec:linear-model:mean-overlap}, and Jensen's inequality, we obtain
\begin{align*}
   \ba_n(0,x_0,\l) - \ba_n(0,x_0,0)&=\unsur{2n} \gesp{\log \crochet{e^{\l n (q-\qud)^2}}}\\
&\le \unsur{2n}\log \gesp{\crochet{e^{\l n (q-\qud)^2}}} \le \unsur{2n} \log L\,.
\end{align*}

}

\subsubsection*{Fact 2}

We shall now use the following fact, coming from the proof of Theorem~\ref{thm:sk:bornesup} (see equation~(\ref{eq:sk:calculderivee}): for every $x_0\ge 2 qt$,

\begin{equation}
  \label{eq:faittrois}
  \ta_n(t,x_0-2qt) = \ta_n(0,x_0) + \frac t2 q^2 -\undemi \intot \gesp{\crochet{(q-\qud)^2}_{s,x_0-2qs,0}}\, ds\,.
\end{equation}

We can now introduce the function 
$$  h_n(t) \egaldef \ta_n(0,x_0) + \frac t2 q^2 -\ta_n(t,x_0-2qt)=\intot \gesp{\crochet{(q-\qud)^2}_{s,x_0-2qs,0}}\, ds\,.$$

For $\l,t,x_0$ small enough, and $q=\qlin(x_0)$, we have

\begin{align*}
\l \frac{d}{dt} h_n(t) &= \frac{\l}{2}\gesp{\crochet{(q-\qud)^2}_{t,x_0-2qt,0}}&\text{(from (\ref{eq:faittrois}))}\\
&\le \unsur{2n} \gesp{\log \crochet{e^{\l n(q-\qud)^2}}_{t,x_0-2qt,0}}&\text{(Jensen's inequality)}\\
&=  \ba_n(t,x_0-2qt,\l) - \ba_n(t,x_0-2qt,0) \\
&\le \ba_n(0,x_0,\l+t) + \frac{t}{2} q^2 - \ba_n(t,x_0-2qt,0)&\text{(from (\ref{eq:faitun}))}\\
&\le \ba_n(0,x_0,0) + \unsur{2n} \log L  \frac{t}{2} q^2 - \ba_n(t,x_0-2qt,0)&\text{(from Proposition~\ref{sec:linear-model:mean-overlap})}\\
&= h_n(t) +\unsur{n} L\,.
\end{align*}
From this differential inequality, we easily get that for $t,x_0,\l$ small enough, and $q=\qlin(x_0)$, we have
$$ h_n(t) \le e^{t/\l} \unsur{n} L \,.$$

Since $h_n$ is positive and $h_n(0)=0$, by taking  $\l=t$ we obtain that $h_n(t) \to 0$.

 We now impose the constraint, $x_0=2qt$ and that is possible because $q=q_c$ is the solution of the equation $q_c = \qlin(2q_c t)$, and then $h_n(t) \to 0$ means exactly for $q=q_c(t)$
$$\ta_n(0,2qt) + \frac t2 q^2 -\ta_n(t,0)= \alphalin(2qt) + \frac t2 q^2 -\alpha_n(t) \to 0\,.$$
\end{proof}
\section{A framework for Gaussian integration by parts}
\label{sec:ipp}
Let $X$ be a random variable defined on the probability space $(\Omega,\Frond,\PP)$ with values in the configuration space $(\Gamma,\Frond^\Gamma)$.

We assume given a gaussian environment, that is another probability space 
 $(\Omega^{(g)},\Frond^{(g)},\gPP)$ and a bimeasurable Hamiltonian
$ H: \Omega^{(g)} \time \Gamma \to \R$ such that $(H(\gamma),\gamma\in\Gamma)$ is a centered Gaussian process with covariance 
$$\gesp{H(\gamma_1)H(\gamma_2)} = R(\gamma_1,\gamma_2)\,.$$

For example, in the classical \SK model, $\Gamma=\ens{-1,1}^n$ is the space of configurations of n spins $\sigma=(\sigma_1, ... ,\sigma_n)$ with $\sigma_i =\pm 1$ and under $\PP$ the spins are independent Bernoulli ($\pm 1$) random variables, and $H_n(\sigma) = \sqrt{\frac{2}{n}} \sum_{i<j} \sigma_i \sigma_j g_{ij}$ with $g_{ij}$ iid standard Gaussian : 
$$ R(\sigma,\tau) = n(\qud^2(\sigma,\tau) -\unsur{n} \qud(\sigma^2,\tau^2))\,.$$

To use the replica technique, we assume given iid configurations $(X_i)_{i\ge 1}$ defined on the same probability space. We can then consider the Gibbs measure of one or several independent replicas, with respect to the same environment:
\begin{gather*}
\crochet{\phi(X)} = \unsur{Z(t)}\esp{\phi(X) e^{\sqrt{t} H(X) -\frac{t}2 R(X,X)}} \\
\crochet{\psi(X_1,X_2)} = \unsur{Z(t)^2}\esp{\psi(X_1,X_2) e^{\sqrt{t} (H(X_1)+H(X_2)) -\frac{t}2 (R(X_1,X_1)+R(X_2,X_2))}} 
\end{gather*}
where $Z(t)$ is the partition function
$$ Z(t) =\esp{ e^{\sqrt{t} H(X) -\frac{t}2 R(X,X)}}\,.$$

\begin{proposition}
If $\alpha(t) = \gesp{\log Z(t)}$ then
$$ \frac{d\alpha}{dt}= -\undemi \gesp{\crochet{R(X_1,X_2)}}\,.$$
\end{proposition}
\begin{proof}
We recall the integration by parts formula (see e.g. Talagrand~\cite{talsaintflour2000}):
if $g$ is centered normal, if the function $f$ is $C^1$, and for some constant $C$, $\valabs{f(x)} \le e^{C\valabs{x}}$, then

$$\esp{g f(g)} = \esp{g^2} \esp{f'(g)}\,.$$
This can be easily extended to functions of several variables. Let $F:\R^n \to \R$ be $C^1$ and such that for a constant $C$ : $\valabs{F(x)} \le e^{C\valabs{x}}$. Then if $(u,u_1, \ldots, u_n)$ is a centered Gaussian vector:

$$\esp{u F(u_1,\ldots,u_n)} = \sum_{i=1}^n \esp{u u_i}\,\esp{\frac{\partial F}{\partial x_i}(u_1, \ldots,u_n)}\,.$$

Differentiating with respect to $t$ yields
$$\frac{d\alpha}{dt} = \unsur{2 \sqrt{t}} \gesp{\crochet{H(X)}} - \undemi \gesp{\crochet{R(X,X)}}\,.$$
Since,
$$ \crochet{H(X)} = \int \prob{X\in dx} \gesp{ H(x) e^{\sqrt{t} H(x) -\frac{t}{2} R(x,x)}\unsur{Z(t)}} = \int \prob{X\in dx} \gesp{ H(x) F(H(\gamma), \gamma \in \Gamma)}$$

for a $C^1$ function $F$ with at most exponential growth, we have
\begin{align*}
&\gesp{\crochet{H(X)}}=  \int \prob{X_1\in dx_1} R(x_1,x_1) \sqrt{t} \gesp{e^{\sqrt{t} H(x_1) - \frac{t}{2}R(x_1,x_1)} \unsur{Z(t)}}\\
& -\int \prob{X_1\in dx_1}\int \prob{X_2 \in dx_2}\sqrt{t} R(x_1,x_2)
\gesp{e^{\sqrt{t} (H(x_1)+H(x_2)) - \frac{t}{2}(R(x_1,x_1)+R(x_2,x_2))}\unsur{Z(t)^2}}\\
& \qquad \ = \sqrt{t}\gesp{\crochet{R(X,X)} - \crochet{R(X_1,X_2)}}\,.
\end{align*}

\end{proof}
With the same type of computations, we obtain
\begin{proposition}
Let $$\tilde{Z}(t) = \esp{e^{\sqrt{t} (H(X_1)+H(X_2)) -\frac{t}2 (R(X_1,X_1)+R(X_2,X_2))+c(X_1,X_2)}} $$ where $c(x_1,x_2)$ is a coupling function,
and let $\tilde{\alpha}(t) = \gesp{\log \tilde{Z}(t)}$. Then
$$\frac{d\tilde{\alpha}}{dt} = \gesp{\crochet{R(X_1,X_2) -2 R(X_1,X_3)}}\,,$$
where in this expression $X_1,\ldots,X_4$ are independent copies of $X$, considered under the same environment, with couplings between $X_1,X_2$ and $X_3,X_4$.
\end{proposition}
\begin{proof}
It is mutatis mutandis the same proof, this time for the Gaussian process
$ K(\gamma_1,\gamma_2) = H(\gamma_1) + H(\gamma_2)$. The expression of the derivative uses symmetry between the $X_i$'s.
\end{proof}

\bibliographystyle{siam}
\bibliography{new,local}

\end{document}